\documentclass[onecolumn,leqno]{amsproc}%
\usepackage{amsfonts}
\usepackage{amsmath}
\usepackage{amssymb}
\usepackage{graphicx}
\usepackage[none]{hyphenat}%
\theoremstyle{plain}

\numberwithin{equation}{section}
\allowdisplaybreaks[3]
\begin{document}
\title{}

\begin{center}
{\LARGE \textbf{Least-Squares Prices of Games}}

\bigskip

\textbf{Yukio Hirashita}

\bigskip

\textbf{Abstract}
\end{center}

\noindent What are the prices of random variables? In this paper, we define
the least-squares prices of coin-flipping games, which are proved to be
minimal, positive linear, and arbitrage-free. These prices depend both on a
set of games that are available for investing simultaneously and on a
risk-free interest rate. In addition, we show a case in which the
mean-variance portfolio theory is inappropriate in our incomplete market.



\bigskip

\begin{center}
{\large \textbf{1. Introduction}}
\end{center}

\noindent Consider the following two coin-flipping games:\smallskip

\noindent\textbf{Game }$\mathbf{A:}$ Profit is $19$ or $1$ if a tossed coin
yields heads or tails, respectively.\smallskip

\noindent\textbf{Game }$\mathbf{B:}$ Profit is $10$ if a tossed coin yields
heads or tails.\smallskip

In general, game $B$ is preferable to game $A$ (see [7\textbf{,} Example
9.2]). Despite the fact that the expectations concerning the two games are
equal, the price of $B$ should be higher than that of $A$. However, if game
$C$ is available for investing simultaneously, the three prices of these games
should be the same; this is because the mixed game $(A+C)/2$ is equal to
$B.$\smallskip

\noindent\textbf{Game }$\mathbf{C:}$ Profit is $1$ or $19$ if a tossed coin
yields heads or tails, respectively.\smallskip

Therefore, the price of a game should change in accordance with the set of
games that are available for investing simultaneously. As F. Black and M.
Scholes demonstrate, the price of an option depends on the risk-free
continuously compound interest rate $r>0$ (see [1,\textbf{ }page 643]). It is
noteworthy that if $r=0$, no investor will invest his/her money, because no
gain is expected. In this paper (except in Remarks 3.5 and 3.6), we assume
that $r$ is $0.05$. The term \textquotedblleft
arbitrage-free\textquotedblright\ implies that no investor has an opportunity
to earn a profit exceeding the risk-free interest rate.

Here, we introduce the pricing method of a coin-flipping game.

\bigskip

\noindent\textbf{Theorem 1.1.} \textit{Suppose that a game }$A:=(a,$\textit{
}$b)$\textit{ involves a profit }$a$\textit{ or }$b$\textit{ }$(a,$\textit{
}$b>0)$\textit{ if a tossed coin yields heads or tails, respectively. Put
}$E^{A}:=(a+b)/2$\textit{. If }$E^{A}/\sqrt{ab}\leq e^{r}$\textit{, then the
price of game }$A$\textit{ is given by }$u_{r}^{A}=\sqrt{ab}/e^{r}$\textit{,
and the optimal proportion of investment is }$1$\textit{. Otherwise, }%
$u_{r}^{A}=\kappa a+(1-\kappa)b$ $($\textit{if }$a\geq b)$\textit{ or }$\kappa
b+(1-\kappa)a$ $($\textit{if }$a<b)$\textit{, where }$\kappa$\textit{
}$:=(1-\sqrt{1-1/e^{2r}})/2$\textit{, and the optimal proportion of investment
is }$u_{r}^{A}(E^{A}$\textit{ }$-u_{r}^{A})$\textit{ }$/((a-u_{r}^{A}%
)$\textit{ }$(u_{r}^{A}-b))$\textit{.}

\noindent\textbf{Proof}\textit{.} Using Remark 3.1 under the conditions of
this theorem and solving the simultaneous quadratic equations, we obtain the
conclusion.\hfill$\square$\smallskip

\noindent\textbf{Game }$\mathbf{A:}$ As $E^{A}/\sqrt{ab}=10/\sqrt
{19}\fallingdotseq2.294>e^{0.05}\fallingdotseq1.051$ and $\kappa
\fallingdotseq0.3458$, we obtain the price $u_{r}^{A}\fallingdotseq7.224$ and
the optimal proportion of investment $t_{u_{r}^{A}}\fallingdotseq
0.274$.\smallskip

Now, we explain the term \textquotedblleft optimal proportion of
investment.\textquotedblright\ Let $t\in\lbrack0,$ $1]$ be a proportion of
investment; then, the investor repeatedly invests $t$ of his/her current
capital. For example, let $c$ be the current capital; when the investor plays
game $A=(19,$ $1)$ once, his/her capital will be $19ct/u+c(1-t)$ or
$ct/u+c(1-t)$ if a tossed coin yields heads or tails, respectively, where $u$
is the price of the game. Let the initial capital be $1$. After $N$ attempts,
if the investor has capital $c_{N}$, then the growth rate (geometric mean) is
given by $c_{N}^{1/N}$. As the value $\lim_{N\rightarrow\infty}\left(
\text{expectation of }c_{N}^{1/N}\right)  $ is a function with respect to
$t\in\lbrack0,$ $1]$, it reaches its maximum of $G_{u}$ at $t=t_{u}$, where
$G_{u}$ is called the limit expectation of the growth rate. The price
$u_{r}^{A}$ is determined by the equation $G_{u}=e^{r}$. Thus, the optimal
proportion of investment $t_{u_{r}^{A}}$ is determined (see Remark 3.1). It is
noteworthy that the value $\lim_{N\rightarrow\infty}\left(  \text{variance of
}c_{N}^{1/N}\right)  $ is $0$.\smallskip

\noindent\textbf{Game }$\mathbf{B:}$ As $E^{B}/\sqrt{ab}=10/\sqrt
{100}=1<e^{0.05}\fallingdotseq1.051$, we obtain the price $u_{r}^{B}$
$\fallingdotseq9.512$. In this case, the optimal proportion of investment is
$1$. This implies that the investor should invest his/her entire current
capital in each attempt.\smallskip

In Section 2, we will introduce the least-squares price $u_{r}^{A\text{,
}\Omega}$ of game $A$ in a set $\Omega$ of games and prove some properties of
$u_{r}^{A\text{, }\Omega}$.\medskip

\noindent\textbf{Example 1.1.} $\Omega=\{(19,1),$ $(4,16)\}$. Using Theorem
1.1, we have $u_{r}^{(19,\text{ }1)}\fallingdotseq7.224$ and $u_{r}^{(4,\text{
}16)}\fallingdotseq8.149$. As $0.4(19,$ $1)+0.6(4,16)=(10,10)$, using Lemma
2.2, we obtain $u_{r}^{(19,\text{ }1),\text{ }\Omega}$ $=u_{r}^{(4,\text{
}16),\text{ }\Omega}=10/e^{r}\fallingdotseq9.512$, where each price reaches
its maximum.\medskip

\noindent\textbf{Example 1.2.} $\Omega=\{(19,$ $1),$ $(16,$ $4)\}$. As Example
1.1, we have $u_{r}^{(19,\text{ }1)}\fallingdotseq7.224$ and $u_{r}%
^{(16,\text{ }4)}\fallingdotseq8.149$. Observe that $p(19,$ $1)+(1-p)(16,$
$4)=(3p+16,$ $4-3p)$ and $10/\sqrt{(3p+16)(4-3p)}\geq5/4>e^{0.05}%
\fallingdotseq1.051$ for each $p\in\lbrack0,1]$. In this case, using Lemma 2.3
with the linearity $u_{r}^{(3p+16,\text{ }4-3p)}$ $=6p\kappa+12\kappa$
$-3p+4$, we obtain $u_{r}^{(19,\text{ }1),\text{ }\Omega}$ $=u_{r}^{(19,\text{
}1)}$ and $u_{r}^{(16,\text{ }4),\text{ }\Omega}=u_{r}^{(16,\text{ }4)}$,
where each price is unchanged.\medskip

\noindent\textbf{Example 1.3.} $\Omega=\{(12,$ $8),$ $(11,$ $9)\}$. Using
Theorem 1.1, we have $u_{r}^{(12,\text{ }8)}\fallingdotseq9.320$ and
$u_{r}^{(11,\text{ }9)}\fallingdotseq9.465$. Observe that $p(12,$
$8)+(1-p)(11,$ $9)=(p+11,$ $9-p)$ and $10/\sqrt{(p+11)(9-p)}\leq5/(2\sqrt{6})$
$\fallingdotseq1.021$ $<e^{0.05}\fallingdotseq1.051$ for each $p$ $\in
\lbrack0,1]$. In this case, by the fact that $u_{r}^{(p+11,\text{ }9-p)}%
=\sqrt{(p+11)(9-p)}/e^{r}$ and using numerical calculations according to
Definition 2.1, we obtain $u_{r}^{(12,\text{ }8),\text{ }\Omega}%
\fallingdotseq9.345$ and $u_{r}^{(11,\text{ }9),\text{ }\Omega}\fallingdotseq
9.469$, where $u_{r}^{A}<u_{r}^{A,\text{ }\Omega}$ $<E^{A}/e^{r}$ for each
$A\in\Omega$.

\bigskip

For a better understanding of the background, we present our incomplete market
assumptions as follows (compare with [8, Sections 2.1 and 8.2]).

1. Frictionless Market: There are no transactions costs or taxes, and all
securities are perfectly divisible.

2. Price-Taker: The investor's actions cannot affect the probability
distribution of returns on the securities. Every security has a positive expectation.

3. No Arbitrage Opportunities: There exits a unique riskless standard asset,
that is not necessarily tradable. Further, there exists a security, wherein
the limit expectation of the growth rate is equal to that of the riskless
standard asset. The limit expectation of the growth rate of any security never
exceeds that of the standard asset. The standard asset is usually provided by
the riskless rate of interest.

4. No Short Sales: Combined with suitable transactions, all necessary short
sales must be included in the securities (probability distribution of
returns), that have positive expectations. For example, $-(19,$ $1)+2(16,$
$4)=(13,$ $7)$.

\bigskip

\begin{center}
{\large \textbf{2. Least-Squares Prices}}
\end{center}

\noindent Let $\Psi:=\{G_{j}:=(c_{j},$ $d_{j})$ : $c_{j},$ $d_{j}>0,$
$j=1,2,...,m\}$ be a finite set of coin-flipping games, which are completely
correlated. Denote the convex cone $\{\sum_{j=1}^{m}k_{j}G_{j}$ : $k_{j}%
\geq0,$ $j=1,2,...,m\}$ by $\widehat{\Psi}$. Then, a basis $\Omega:=$
$\{A_{i}$ : $i=1$ or $i=1$, $2\}$ exists such that $\widehat{\Psi}%
=\widehat{\Omega}$ (see Remark 3.3). This is because, if $\min_{j=1,2,...,m}%
c_{j}/d_{j}$ $=\max_{j=1,2,...,m}c_{j}/d_{j}$, we can choose $\Omega
:=\{A_{1}:=G_{1}\}$. If not, we can choose $\Omega:=$ $\{A_{1}:=G_{j_{0}},$
$A_{2}:=G_{j_{1}}\}$ such that $c_{j_{0}}/d_{j_{0}}$ $=\min_{j=1,2,...,m}$
$c_{j}/d_{j}$ and $c_{j_{1}}/d_{j_{1}}$ $=\max_{j=1,2,...,m}$ $c_{j}/d_{j}$.
Since the set $\Omega=\{A_{i},$ $1\leq i\leq n\}$ $(n=1$ or $2)$ is a basis of
the convex cone $\widehat{\Omega}$, if $\sum_{i=1}^{n}$ $k_{i}A_{i}=\sum
_{i=1}^{n}k_{i}^{\prime}A_{i}$, then $k_{i}=k_{i}^{\prime}$ for each $1\leq
i\leq n$.

Set $S:=\{(t_{i})\in R^{n}$ : $0\leq t_{i}\leq1$, $1\leq i\leq n\}$ and $Q:=$
$\{(p_{i})\in R^{n}$ : $\sum_{i=1}^{n}p_{i}=1$, $p_{i}\geq0$, $1\leq i\leq
n\}$.

From Theorem 1.1, we can verify that $0<u_{r}^{A}\leq E^{A}/e^{r}$ and
$u_{r}^{kA}:=u_{r}^{(ka,\text{ }kb)}$ $=ku_{r}^{A}$ for each $A=(a,$ $b)$ and
$k>0$.\medskip

\noindent\textbf{Definition 2.1.} By defining the function
\begin{equation}
L((t_{i})):=\sup_{(p_{i})\in Q}\frac{u_{r}^{\sum_{i=1}^{n}p_{i}A_{i}}}%
{\sum_{i=1}^{n}p_{i}(u_{r}^{A_{i}}+t_{i}(E^{A_{i}}/e^{r}-u_{r}^{A_{i}}%
))}\text{ \ \ }((t_{i})\in S), \tag{2.1}%
\end{equation}
we have $L((0))=\sup_{(p_{i})\in Q}(u_{r}^{\sum_{i=1}^{n}p_{i}A_{i}}%
/\sum_{i=1}^{n}p_{i}u_{r}^{A_{i}})\geq u_{r}^{A_{1}}/u_{r}^{A_{1}}=1$ and%
\[
L((1))=\sup_{(p_{i})\in Q}\frac{u_{r}^{\sum_{i=1}^{n}p_{i}A_{i}}}{\sum
_{i=1}^{n}p_{i}E^{A_{i}}/e^{r}}=\sup_{(p_{i})\in Q}\frac{u_{r}^{\sum_{i=1}%
^{n}p_{i}A_{i}}}{E^{\sum_{i=1}^{n}p_{i}A_{i}}/e^{r}}\leq1.
\]
Since the set $T:=\{(t_{i})\in S$ : $L((t_{i}))\leq1\}$ is not null, convex,
closed, and thus compact, there is a unique point $(x_{i})\in T$ such that
$v:=\min_{(t_{i})\in T}$ $\sum_{i=1}^{n}$ $t_{i}^{2}=\sum_{i=1}^{n}x_{i}^{2}$.
Define $u_{r}^{A_{i},\text{ }\Omega}:=u_{r}^{A_{i}}$ $+x_{i}(E^{A_{i}}%
/e^{r}-u_{r}^{A_{i}})$ and call it the \textit{least-squares price of }$A_{i}%
$\textit{ in }$\Omega$ for each $1\leq i\leq n.$ For each mixed game
$\sum_{i=1}^{n}k_{i}A_{i}\in\widehat{\Omega}$, define $u_{r}^{\sum_{i=1}%
^{n}k_{i}A_{i},\text{ }\Omega}$ $:=\sum_{i=1}^{n}k_{i}u_{r}^{A_{i},\text{
}\Omega}$ (see [7, Section 9.6]).\medskip

\noindent\textbf{Lemma 2.2.} \textit{If }$\left(  p_{i}\right)  \in Q$\textit{
exists such that }$\sum_{i=1}^{n}p_{i}A_{i}$\textit{ is constant and }%
$p_{k}\neq0$\textit{, then }$u_{r}^{A_{k},\text{ }\Omega}=E^{A_{k}}/e^{r}%
$\textit{.}

\noindent\textbf{Proof}. Write $B:=\sum_{i=1}^{n}p_{i}A_{i}$; then as $B$ is
constant, $u_{r}^{B}=E^{B}/e^{r}$. From Definition 2.1, we obtain
$u_{r}^{A_{i}}\leq u_{r}^{A_{i},\text{ }\Omega}\leq E^{A_{i}}/e^{r}$ and
$u_{r}^{B}\leq\sum_{i=1}^{n}p_{i}u_{r}^{A_{i},\text{ }\Omega}\leq\sum
_{i=1}^{n}p_{i}E^{A_{i}}/e^{r}$ $=u_{r}^{B}.$ Thus, $u_{r}^{A_{k},\text{
}\Omega}=E^{A_{k}}/e^{r}$ if $p_{k}\neq0$.$\hfill\square$\medskip

\noindent\textbf{Lemma 2.3. }\textit{If }$u_{r}^{\sum_{i=1}^{n}p_{i}A_{i}%
}=\sum_{i=1}^{n}p_{i}u_{r}^{A_{i}}$\textit{ for each }$(p_{i})\in Q$\textit{,
then }$u_{r}^{A_{i},\text{ }\Omega}$\textit{ }$=u_{r}^{A_{i}}$\textit{
}$(1\leq i$\textit{ }$\leq n)$\textit{.}

\noindent\textbf{Proof}. By the above assumption, we obtain $L((0))=1$ and
$v=0$, which implies the conclusion.$\hfill\square$\medskip

\noindent\textbf{Theorem 2.4. }\textit{The system of least-squares prices is
arbitrage-free, and there is a mixed game that earns profit equal to the
growth rate of }$e^{r}$\textit{.}

\noindent\textbf{Proof}. As $T\subset S$ and $Q$ are compact, and $u_{r}%
^{\sum_{i=1}^{n}p_{i}A_{i}}$ is continuous with respect to $(p_{i})\in Q$ (see
Theorem 1.1), $(x_{i})\in T$ and $(q_{i})\in Q$ exist such that
\begin{equation}
L((x_{i}))=\max_{(p_{i})\in Q}\frac{u_{r}^{\sum_{i=1}^{n}p_{i}A_{i}}}%
{\sum_{i=1}^{n}p_{i}(u_{r}^{A_{i}}+x_{i}(E^{A_{i}}/e^{r}-u_{r}^{A_{i}}%
))}=\frac{u_{r}^{\sum_{i=1}^{n}q_{i}A_{i}}}{\sum_{i=1}^{n}q_{i}u_{r}%
^{A_{i},\text{ }\Omega}}=1. \tag{2.2}%
\end{equation}
This shows that the mixed game $\sum_{i=1}^{n}q_{i}A_{i}$ earns profit that is
equal to the growth rate of $e^{r}$. On the other hand, for each nonzero mixed
game $\sum_{i=1}^{n}k_{i}A_{i}$ $=k\sum_{i=1}^{n}p_{i}A_{i}$ $\in
\widehat{\Omega}$ $(k:=\sum_{i=1}^{n}k_{i}$, $k>0,$ $p_{i}:=k_{i}/k$,
$(p_{i})\in Q)$, by equation (2.2), we have $u_{r}^{\sum_{i=1}^{n}k_{i}A_{i}}$
$=ku_{r}^{\sum_{i=1}^{n}p_{i}A_{i}}\leq k\sum_{i=1}^{n}p_{i}u_{r}%
^{A_{i},\text{ }\Omega}=\sum_{i=1}^{n}k_{i}u_{r}^{A_{i},\text{ }\Omega}$.
Therefore, the game $\sum_{i=1}^{n}k_{i}A_{i}$ earns profit that is equal to
or less than the growth rate of $e^{r}$.\hfill$\square$\medskip

\noindent\textbf{Theorem 2.5. }\textit{The system of least-squares prices is
minimal in order to be arbitrage-free}.

\noindent\textbf{Proof.} We prove this by using reduction to absurdity.
Assuming that a set of prices $\{R_{i}\}$ of $\{A_{i}\}$ exists such that
$R_{i}\leq u_{r}^{A_{i},\text{ }\Omega}$ $(1\leq i\leq n)$ and $R_{k}$
$<u_{r}^{A_{k},\text{ }\Omega}$ for some $k$. If $R_{j}<u_{r}^{A_{j}}$ for
some $j$, then the game $A_{j}$ earns profit exceeding the growth rate of
$e^{r}$. Thus, we can assume that $u_{r}^{A_{i}}\leq R_{i}$ $(1\leq i$ $\leq
n)$. Therefore, $(s_{i})\in S$ exists such that $R_{i}=u_{r}^{A_{i}}$
$+s_{i}(E^{A_{i}}/e^{r}-u_{r}^{A_{i}})$, where $s_{j}:=0$ is chosen if
$u_{r}^{A_{j}}=E^{A_{j}}/e^{r}$. It is easy to verify that $s_{i}\leq x_{i}$
$(1\leq i\leq n)$ and $s_{k}<x_{k}$. From the above statement, we have
$\sum_{i=1}^{n}$ $s_{i}^{2}$ $<\sum_{i=1}^{n}x_{i}^{2}$, which implies
$(s_{i})$ $\notin T$, and thus $L((s_{i}))>1$. Therefore, a point $(q_{i})$
$\in Q$ exists such that $\sum_{i=1}^{n}$ $q_{i}R_{i}$ $<u_{r}^{\sum_{i=1}%
^{n}q_{i}A_{i}}$, that is, the mixed game $\sum_{i=1}^{n}q_{i}A_{i}$ earns
profit exceeding the growth rate of $e^{r}$.\hfill$\square$\medskip

It is not difficult to verify that if $\Omega=\{A_{i}$ : $1\leq i\leq n\}$ and
$\Omega^{^{\prime}}=\{B_{j}$ : $1\leq i\leq s\}$ are the bases of the convex
cone $\widehat{\Psi}$, then $n=s$, $A_{i}=v_{i}B_{i}$, and $u_{r}%
^{A_{i},\text{ }\Omega}=v_{i}u_{r}^{B_{i},\text{ }\Omega^{\prime}}$ $(v_{i}%
>0$, $1\leq i$ $\leq n)$ after the permutations. Therefore, $u_{r}^{A,\text{
}\Omega}=u_{r}^{A,\text{ }\Omega^{\prime}}$ for each $A\in\widehat{\Psi}$, and
we can define $u_{r}^{A,\text{ }\Psi}:=u_{r}^{A,\text{ }\Omega}$.

\bigskip

\begin{center}
{\large \textbf{3. Remarks}}
\end{center}

\bigskip

\noindent\textbf{Remark 3.1.} Consider a random variable $X$ with nonnegative
bounded profit $a(x)$ and distribution $dF(x)$. In the case where $\exp
(\int\log a(x)dF(x))/e^{r}$ $\leq1/\int1/a(x)dF(x)$, the price is given by
$u_{r}^{X}=\exp(\int\log a(x)dF(x))/e^{r}$, and the optimal proportion of
investment is $1$. Otherwise, the price $u=u_{r}^{X}$ and the optimal
proportion of investment $t$ are determined by the simultaneous equations
$\exp(\int\log(a(x)t/u-t+1)$ $dF(x))=e^{r}$ and $\int(a(x)-u)/(a(x)t$
$-ut+u)dF(x)=0$ (see [3, Corollaries 5.1, 5.3, and Section 6]).\medskip

\bigskip

\noindent\textbf{Remark 3.2.} Remark 3.1 can be generalized to the nonnegative
unbounded case where $\int_{a(x)>1}a(x)^{\nu}dF<\infty$ for some $\nu>0$. For
example, because $\sum_{j=1}^{\infty}(2^{j})^{1/2}/2^{j}=2/(2-\sqrt{2})$ and
$\exp(\sum_{j=1}^{\infty}(\log2^{j})/2^{j})/e^{r}=4/e^{0.05}$ $>1/\sum
_{j=1}^{\infty}1/4^{j}=3$, the St. Petersburg game $\{$profit $2^{j}$ with
probability $1/2^{j},$ $\ j=1,2,...\}$ is priced at $4.816$ with the optimal
proportion of investment $0.204.$\medskip

\bigskip

\noindent\textbf{Remark 3.3.} In Section 2, the value of $n$ is $1$ or $2$.
However, when the reader challenges to study dice games, the value of $n$ may
be $36$. To generalize this theory to the convex cone $\widehat{\Psi}$ with a
finite basis $\Omega$, we need the fact that $u_{r}^{\sum_{i=1}^{n}p_{i}A_{i}%
}$ is concave and continuous with respect to $(p_{i})\in Q$ for any positive
integer $n$. This can be achieved using [2, Theorems 185 and 214] and [9,
Theorems 10.1, 10.3 and 20.5] with tedious discussions. Therefore, in
Definition 2.1, $\sup_{(p_{i})\in Q}$ can be replaced by $\max_{(p_{i})\in Q}$
because of Berge's maximum theorem (see [10, Theorem 2.1]).

\bigskip

\noindent\textbf{Remark 3.4. }Let $S$ denote the stock price which is a
nonnegative random variable. Define $P:=\max(K-S,$ $0)$ and $C:=\max(S-K,$
$0)$ for the strike price $K$. Applying Lemma 2.2 with $\Omega=\{P,$ $C,$
$S-C\}$, the equalities $P+(S-C)=K$ and $C+(S-C)=S$ imply Put-call parity
$u_{r}^{C,\text{ }\Omega}-u_{r}^{P,\text{ }\Omega}+K/e^{r}=u_{r}^{S,\text{
}\Omega}$ (see [7, Sections 12.3 and 13.2]).\medskip

\bigskip

\noindent\textbf{Remark 3.5. }In this remark, we assume that the risk-free
interest rate $r=0.02$ is simple (not continuously compound). Consider two
independent coin-flipping games, $X=(50,$ $1)$ and $Y=(30.6191,$ $14)$, where
the variances are $v_{X}=600.25$ and $v_{Y}$ $=$ $69.0486,$ respectively.
Assume that the rates of mean return ([7, section 6.4]) are $r_{X}$
$=0.233546$ and $r_{Y}=0.079211$, respectively. Thus, from the one-fund
theorem ([7, section 6.9]), we have the weight%
\[
w_{X}=\frac{\frac{r_{X}-r}{v_{X}}}{\frac{r_{X}-r}{v_{X}}+\frac{r_{Y}-r}{v_{Y}%
}}=0.2932,
\]
which implies that the single fund of risky assets is
\[
w_{X}X+(1-w_{X})Y=(36.3016,\text{ }24.5552,\text{ }21.9348,\text{ }10.1884).
\]
The four values of this fund occur with the same probability of $1/4$. The
price of this fund is $21.3995$ according to Remark 3.1, where $e^{r}$ is
replaced by $1+r$.

However, the price of $wX+(1-w)Y$ $(0\leq w\leq1)$ reaches its maximum value
of $21.4134$ when $w=0.3514$, that is, the fund
\[
0.3514X+0.6486Y=(37.4295,\text{ }26.6504,\text{ }20.2109,\text{ }9.4318)
\]
is more valuable than the single fund $w_{X}X+(1-w_{X})Y$ because $21.3995$
$\ <$ $\ 21.4134.$

It should be noted that by using $1+r$ instead of $e^{r}$, Theorem 1.1 gives
the prices of $X$ and $Y$ as $u_{X}=20.6721$ and $u_{Y}=20.6721$,
respectively. Thus, the corresponding rates of mean return are $r_{X}$
$=25.5/u_{X}-1$ $=0.233546$ or $r_{Y}=22.30955/u_{Y}-1$ $=0.07921$.

Moreover, Remark 3.1 gives us the optimal proportion $t=0.4222$ for the risky
fund. Thus, the best proportions of investment to $X$, to $Y$, and the
risk-free asset are $tw=0.1484,$ $t(1-w)=0.2738,$ and $1-t=0.5778,$
respectively. The mean-variance portfolio theory cannot provide a proportion
of $0.5778$ for the risk-free asset (see [7, section 7.1]).

\bigskip

\noindent\textbf{Remark 3.6. }Let $Y:=\{Y_{t}\}_{0\leq t\leq T<\infty}$ be a
measurable stochastic process with a filtration. Put $\Psi:=\{\tau$ $;$ $\tau$
is a stopping time such that $\tau\leq T$\}. We define the price $\overline
{u}^{Y}$ by $\sup_{\tau\in\Psi}u_{rE(\tau)}^{Y_{\tau}},$ where $u_{rE(\tau
)}^{Y_{\tau}}$ is the price of the random variable $Y_{\tau(\omega)}(\omega)$
(see [Karatzas et al. (1998)]) with respect to the growth rate $e^{rE(\tau)}$:%
\[
\sup_{_{\substack{0\leq z\leq1\\z\leq\text{ess }\inf_{\omega}\text{ }%
Y_{\tau(\omega)}(\omega)\text{ }z/u+1}}}\int_{\Omega}\log\left(
Y_{\tau(\omega)}(\omega)z/u-z+1\right)  d\omega=rE(\tau).
\]
The geometric price of $Y$ $:=\{Y_{t}\}_{0\leq t<\infty}$ is defined by
$\sup_{0<T<\infty}$ $\overline{u}^{\{Y_{t}\}_{0\leq t\leq T}}$.

\newpage

\begin{center}
\textbf{References}
\end{center}

\noindent\lbrack1] F. Black and M. Scholes, The pricing of options and
corporate liabilities,

\textit{Journal of Political Economy} \textbf{81} (1973) 637--654.

\noindent\lbrack2] G. H. Hardy, J. E. Littlewood and G. P\'{o}lya,
\textit{Inequalities} (Cambridge University

Press, Reprinted 1973).

\noindent\lbrack3] Y. Hirashita, Game pricing and double sequence of random
variables, Preprint,

arXiv:math.OC/0703076 (2007).

\noindent\lbrack4] Y. Hirashita, Delta hedging without the Black-Scholes
formula, \textit{Far East }

\textit{Journal of Applied Mathematics} \textbf{28} (2007) 157-165.

\noindent\lbrack5] I. Karatzas and S. E. Shreve, \textit{Brownian motion and
stochastic calculus},

(Springer-Verlag, New York, 1998).

\noindent\lbrack6] J. L. Kelly, Jr., A new interpretation of information
rate,\ \textit{Bell System Technical }

\textit{Journal} \textbf{35} (1956) 917--926.

\noindent\lbrack7] D. G. Luenberger, \textit{Investment science} (Oxford
University Press, Oxford, 1998).

\noindent\lbrack8] R. C. Merton, Continuous-time finance (Basil Blackwell,
Cambridge, 1990).

\noindent\lbrack9] R. T. Rockafellar, \textit{Convex analysis} (Princeton
University Press, Princeton,

1970).

\noindent\lbrack10] I. E. Schochetman, Pointwise versions of the maximum
theorem with applications

in optimization, \textit{Applied Mathematics Letters} \textbf{3} (1990) 89--92.

\bigskip

Chukyo University, Nagoya 466-8666, Japan

yukioh@lets.chukyo-u.ac.jp

\bigskip

To help the readers to understand this article,

\ \ \ \ Y. Hirashita, Details for \textquotedblright Least-Squares Prices of
Games.\textquotedblright\ Preprint, 2007,

\ \ \ \ arXiv:math.OC/0703080

is posted.
\end{document}